\input amstex
\magnification =1200
\documentstyle {amsppt}

\define\pd{\partial}

\define\spsh{strictly plurisubharmonic }
\define\gt{Grauert tube}
\define\Rm{Riemannian manifold}
\document
\baselineskip 1.4\baselineskip
\topmatter
\title  On rigidity of Grauert tubes over homogeneous Riemannian manifolds\endtitle
\author Su-Jen Kan 
\endauthor 
\address 
Institute of Mathematics, Academia Sinica, Taipei 115, Taiwan
\endaddress
\email 
{kan\@math.sinica.edu.tw}
\endemail
\date Dec. 21, 2003
\enddate
\NoRunningHeads
\thanks 2000 {\it Mathematics Subject Classification. } 32C09, 32Q28, 32Q45, 53C24, 58D19. This
research is partially supported by NSC 90-2115-M-001-006. 
\endthanks
\abstract

Given a real-analytic Riemannian manifold $X$ there is a canonical  complex structure, which
is compatible with the canonical complex structure on $T^*X$ and  makes the leaves of the
Riemannian foliation on
$TX$ into holomorphic curves, on  its tangent bundle. A {\it Grauert tube} over $X$ of radius $r$,
denoted as $T^rX$, is the collection of tangent vectors of $X$ of length less than $r$ equipped
with this canonical complex structure. 

In this article, we prove the following two rigidity property of Grauert tubes. First, for any
real-analytic Riemannian manifold such that $r_{max}>0$, we show that the identity component of
the automorphism group of $T^rX$ is isomorphic to the identity component of the isometry group
of $X$ provided that $r<r_{max}$. Secondly, let $X$ be a
homogeneous Riemannian manifold and let the radius $r<r_{max}$, then the automorphism group of
$T^rX$ is isomorphic to the isometry group of $X$ and there is a unique Grauert tube
representation for such a complex manifold $T^rX$.

\endabstract
\endtopmatter

\subheading {1. Introduction}

The purpose of this article is to give an affirmative answer to the rigidity of \gt s over
(compact or non-compact) homogeneous Riemannian  spaces. On the way to prove this, we are also
able to show that for any real-analytic Riemannian manifold,  the identity component of
the automorphism group of a
\gt\ is isomorphic to the identity component of the isometry group of the center when it is
defined and the radius is not the critical one.

It was  observed by Grauert [6] that  a real-analytic manifold $X$ could be embedded in
 a complex manifold as a maximal totally real submanifold. One way to see this is to
complexify the transition functions defining $X$. However, this complexification is not
unique. In [7] and [12], Guillemin-Stenzel and independently Lempert-Sz\H oke encompass
certain conditions on the ambient complex structure to make the complexification canonical
for a given real-analytic Riemannian manifold. In short, they were looking for a complex
structure,  on part of the cotangent bundle
$T^*X$,  compatible with the canonical symplectic structure on $T^*X$. Equivalently,
it is to say that there is a unique complex structure, called the adapted complex
structure, on part of the tangent bundle of $X$ making the leaves of the Riemannian
foliation on $TX$ into  holomorphic curves. The set of tangent vectors of length less than
$r$ equipped with the adapted complex structure is called a \gt\ $T^rX$. 

Since the
adapted complex structure is constructed canonically associated to the  metric
$g$ of the Riemannian manifold $X$, the differentials of the isometries of $X$ are
automorphisms of $T^rX$. Conversely, it is interesting to see whether all automorphisms of
$T^rX$ come from the differentials of the isometries of $X$ or not. When the answer is
affirmative, we say the \gt\ is {\it rigid}. 

With respect to the adapted complex structure, the length square function $\rho
(x,v)=|v|^2, v\in T_xX$, is strictly plurisubharmonic. When the center $X$ is compact, the
\gt \ $T^rX$ is exhausted by $\rho$, hence is a Stein manifold with smooth strictly
pseudoconvex boundary when the radius is less than the critical one. In this case, the
automorphism group of $T^rX$ is a compact Lie group. Base on these, Burns, Burns-Hind ([1],
[2]) are able to prove the rigidity for \gt s over compact real-analytic Riemannian
manifolds. In [1], Burns also showed that the rigidity of a \gt\ is equivalent to the
uniqueness of center when the \gt\ is constructed  over a compact manifold.  

When $X$ is non-compact nothing particular is known, not even to the general existence of a
\gt\ over $X$. When it exists, most of the good properties in the compact cases are
lacking here since the length square function $\rho$ is no longer an exhaustion.
However, in this article we are able to prove some kind of weak rigidity for general Grauert
tubes. That is, we are able to show that $Aut_0(T^rX)=Isom_0(X)$ when it is defined and
$r<r_{max}$.

 By now, there are very few complete non-compact cases we are sure about the existence of \gt s.
One kind  of them are those over co-compact real-analytic Riemannian manifolds, the \gt s\ are
simply the lifting of the \gt s\ over their  compact quotients. The second kind are Grauert tubes
over real-analytic homogeneous  Riemannian manifolds. In [10], Kan and Ma have proved the
rigidity of
\gt s over (compact or non-compact) symmetric spaces based on the co-compactness of symmetric
spaces. Later on, in [9], the author proves that  the uniqueness of such \gt s\ also holds.

We prove in this article that not only the rigidity but also the uniqueness holds  for \gt s
over real-analytic homogeneous Riemannian manifolds. The center manifold could be either
compact or non-compact. Moreover, it is also true for \gt s constructed over quotient
manifolds of a homogeneous Riemannian space. The main result is Theorem 7.6. The condition
$r<r_{max}$ is necessary since \gt s over non-compact symmetric spaces of rank one of
maximal radius were proved in [3] to be non-rigid.

On the  way to
this we, omitting the requirement on the tautness, prove a generalized version of the
Wong-Rosay theorem on the characterization of ball to any domain in a complex manifold on Theorem
3.1. We also prove at Theorem 5.2  on the complete hyperbolicity of \gt s over
homogeneous Riemannian manifolds. On Theorem 6.4, a weak version of the rigidity is given for
any general Grauert tubes.

The outline of this article is as follows. In $\S 2$,  the  behavior of 
the Kobayashi metric near   $C^2$  strictly pseudoconvex  boundary points was observed. $\S
3$ is devoted to the generalized Wong-Rosay theorem. We give  a brief review on \gt s  and prove one of the key theorem, Theorem 4.1, toward the 
rigidity in $\S 4$. In $\S 5$, the
complete hyperbolicity of \gt s over homogeneous Riemannian  spaces was claimed. $\S 6$ is on
the discussion of the equality of $Aut_0(T^rX)$ and $Isom_0(X)$.   In $\S 7$, we prove the
rigidity
$Isom  (X)=Aut (T^rX)$ and the uniqueness of
\gt s. 
 
The author would like to thank Kang-Tae Kim for pointing out to me that my previous version of 
Rosay-Wong theorem might be able to be generalized to such a generality; to Akio Kodama for
informing me a mistake in the first draft and to Eric Bedford for showing me how to use the
bumping function in the proof of Lemma 2.1. 

\

\subheading {2. Estimates of Kobayashi metrics near a strictly pseudoconvex point}

For complex manifolds  $N$ and $M$, let $Hol(N,M)$ be the family of 
holomorphic mappings from $N$ to $M$. Let $\Delta $ denote the unit
disk in $\Bbb C$.  The {\it Kobayashi  pseudometric } on $M$ is defined by
$$
F_{\sssize M}(z,\xi)
=\inf\{\alpha|\ \alpha>0,\exists f\in H(\Delta,M), f(0)=z,f'(0)=\frac {\xi}{\alpha}\}. $$
The Kobayashi  pseudometric $F_{\sssize M}$ is holomorphic decreasing, i.e., if $f:N\to M$
is a holomorphic map, then $F_{\sssize M}(f(z),f_*\xi)\le F_{\sssize N}(z,\xi).$ We also
quote the following property (c.f. p.128, Proposition 2, [15]) for later use: if $K$ is
a compact subset of $M$ contained in a coordinate polydisk, then there is a constant
$C_{\sssize K}$ such that 
$$F_{\sssize M}(z,\xi)\le C_{\sssize K}\|\xi\|\tag 2.1$$ for all $z\in K$ where $\|\xi\|
=max|\xi_i|$.

Royden [15] has shown that the {\it Kobayashi  pseudodistance} $d_{\sssize
M}^{\sssize K}$ is the integrated form of $F_{\sssize M}$. That is,
given $z,w\in M$,
$$d_{\sssize
M}^{\sssize K}(z, w)=\inf_{\gamma}\int_0^1 F_{\sssize M}(\gamma(t),\gamma'(t))dt$$
where the infimum is taken over all piecewise $C^1$-smooth curves $\gamma:[0,1]\to M$
joining the
 two points $z, w\in M$. 

The manifold $M$ is said to be {\it hyperbolic} if its Kobayashi
pseudodistance is a distance, i.e.,  $d_{\sssize M}^{\sssize K}(z,w)=0$ implies $z=w$.
The hyperbolicity of a manifold $M$ is  equivalent to the following (c.f.[15]): for any $x\in M$
there exists a neighborhood $U$ of $x$ and a positive constant $c$ such that $F_{\sssize
M}(y,\xi)\ge c\|\xi\|$ for all $y\in U$. In this case the Kobayashi pseudodistance is called
the Kobayashi distance of $M$. It is well known that a bounded domain in
$\Bbb C^n$ is hyperbolic. Furthermore, Sibony [16] proved that a complex manifold with a
bounded
\spsh\ function must be hyperbolic. A hyperbolic manifold $M$ is said to be {\it complete
hyperbolic} if its Kobayashi metric is complete, i.e., if each Cauchy sequence in the
Kobayashi metric has a limit or equivalently that any ball of finite radius is relatively
compact in $M$. 

 In [5], Graham  has observed the boundary behavior of 
the Kobayashi metric on  bounded strictly pseudoconvex
domains  $D\subset \Bbb C^n$. He showed that the
boundary is at infinite distance, i.e., $\lim_{q\rightarrow \partial
D}d^{\sssize K}_{\sssize D}(q, K)=\infty$ for any compact $K\subset D$. He hence
concluded that every bounded strictly pseudoconvex domain in $\Bbb C^n$ is complete
hyperbolic.

For the rest of this section, we would like to make a similar estimate for $C^2$-smooth strictly
pseudoconvex boundary points of a domain in a  complex manifold.  
\proclaim{Lemma 2.1}  
Let $D$ be a domain in a  complex manifold $ M$ and let
$q\in\pd D$ be a $C^2$-smooth strictly pseudoconvex
boundary point. Then there exist a domain $D'\supset D$ and 
 a bounded plurisubharmonic function $h$ on $D'$ which is strictly plurisubharmonic on a
neighborhood of $q$.
\endproclaim
\demo{Proof}
Without loss of generality, we may assume that there is a neighborhood $U$ of $q$ in $M$ and a
coordinate chart $\Phi$ from $U$ onto an open set $V\subset \Bbb C^n$ such that $\Phi (q)=0$ and
the local defining function $\rho$ for $D$ at $q$ has the following form:
$$\rho\cdot \Phi ^{-1}(z)= 2 Re\; z_n + |z|^2 + o (|z|^2).$$
Let $$\tilde f(z)=2 Re\; z_n + \frac 12|z|^2 + o (|z|^2)$$ be a plurisubharmonic function  on
$V$.   
 Suitably shrinking the
neighborhood $U$, we see that  the function $f$ defined as 
$f=\tilde f\cdot
\Phi$ is a bounded plurisubharmonic function on $U$ which is strictly plurisubharmonic near the
point $q$. 

Denoting $U^- =U\cap D$, $V^-=\Phi (U^-)$ and $ W^-=\{ w\in U:  f(w)<0\},$ then
 $U^-\subsetneq W^-$ and $f^{-1}(0)\cap \rho^{-1}(0)=\{q\}.$ Therefore,
$$m:= \max _{w\in\overline D\cap (\pd U)} f(w) <0.$$
Let $D'=D\cup U$. We define a function $h$ on $D'$ as following:
$$h=\left \{
\aligned
 \max (f, \frac {m}2)\;\;\;&\text { on }\;\; U,\\
\frac {m}2\;\;\;\;\;\;\;\;\;\;\;\;&\text { on }\; D-U.
\endaligned
\right.$$
It is clear that the domain $D'$ and the function $h$ fulfill the requirement in the statement
of the lemma.\qed
\enddemo

Applying Lemma 2.1 to Theorem 3 of [16], there exist a neighborhood $U$ of $q$ in $D'$ and  a
positive constant $c>0$ such that 
$$F_{\sssize D'}(z,\xi)\ge c\|\xi\|\;\;\text { for all } z\in U, \xi\in
\Bbb C^n.\tag 2.2$$
By the holomorphic decreasing of the Kobayashi pseudodistance, we have 
$$F_{\sssize D}(z,\xi)\ge c\|\xi\|\;\;\text { for all } z\in U\cap D, \xi\in
\Bbb C^n,\tag 2.3$$
which concludes that $M$ is hyperbolic on $U\cap D$.

It was first mentioned by Royden in [15] without proof and was later on proved by Graham at
Lemma 4, p.234 of  [5] that if $P\subset G$ are two complex manifolds (the original requirement
on the hyperbolicity of the manifolds  could be removed away)  then, for
$z\in P$, 
$$F_{\sssize P}(z,\xi)\le \coth (\inf_{w\in G-P}D_{\sssize G}^*(z,w))  F_{\sssize
G}(z,\xi).\tag 2.4$$
 $D_{\sssize G}^*(z,w)$ is defined as follows: 
$$D_{\sssize G}^*(z,w)=\inf\{\kappa (a,b)|\ \exists f\in H(\Delta,G), f(a)=z,f(b)=w\}$$
where $\kappa$ is the Poincar\'e distance for $\Delta$. From the equivalent definition of
the Kobayashi pseudodistance in the classical sense, it is clear that $d_{\sssize G}^{\sssize
K}(z,w)\le D_{\sssize G}^*(z,w).$
 Making a 
little modification on this, we are able to prove the following:
\proclaim{Lemma 2.2}  
Let $D$ be a domain in a  complex manifold $ M$ and let
$q\in\pd D$ be a $C^2$-smooth strictly pseudoconvex
boundary point. If $V$ is a neighborhood of $q$ in $ M$, then
$$\lim_{p\in D, p\to q}d^{\sssize K}_{\sssize D}(p, D\backslash V)=\infty.$$ 
\endproclaim
\demo{Proof}Let the domain $D'$ be constructed as in Lemma 2.1.
By (2.2) there exist a neighborhood $ U$ of $q$ in $ D'$ and a positive constant $c$ such
that $F_{D'}(z,\xi)\ge c\|\xi\|$ for all $z\in  U$. Let $ V\Subset  U$ be
a smaller neighborhood of $q$ in $D'$. Then there exists a $\delta>0$ such that 
$d_{\sssize D'}^{\sssize
K}(z,w)\ge\delta $ for all $z\in  V$ and $w\in  D'- U$. By the holomorphic decreasing property
of the Kobayashi pseudodistance 
$$d_{\sssize D}^{\sssize
K}(z,w)\ge\delta\;\;\;\text { for all }\; z\in  V\cap D,\; w\in  D- U.$$
For all $z\in V\cap D$,
$$\inf_{w\in D-U} D_{\sssize D}^*(z,w)\ge \inf_{w\in D-U}d_{\sssize D}^{\sssize
K} (z,w)\ge\delta>0,$$ and 
$$\coth (\inf_{w\in D-U} D_{\sssize D}^*(z,w))\le coth (\delta)\le L$$
for some constant $L$ independent of $z$.

On the other hand, $U\cap D$ could also
be viewed as a bounded domain $\tilde U$ in $\Bbb C^n$ with $C^2$-smooth strictly
pseudoconvex boundary points through the biholomorphic map $\Phi$. We set $ E=\Phi (V\cap D).$ 
By (2.4), 
$$F_{\sssize D}(z,\xi)\ge \frac 1{L} F_{\sssize U\cap D}(z,\xi)=\frac 1{L}
F_{\sssize \tilde U}(\Phi (z),\Phi_*(\xi))\tag 2.5$$ for all
$z\in V\cap D$.
 The completeness of $d_{\sssize \tilde U}^{\sssize K}$ along with (2.5) leads to,
$$\aligned
\lim_{p\in D, p\to q}d^{\sssize K}_{\sssize D}(p, D\backslash V)&\ge \lim_{p\in
V, p\to q}d^{\sssize K}_{\sssize D}(p, U\backslash V)\\
& \ge \frac 1{L}\lim_{p\in
V, p\to q}d^{\sssize K}_{\sssize \tilde U}(\Phi (p), \tilde U\backslash E)\\&
=\infty.\endaligned$$\qed
\enddemo

\
\

 \subheading {3. An extension of the Wong-Rosay Theorem}

Historically speaking, Wong proved the biholomorphic  equivalence of a bounded
strictly pseudoconvex domain with non-compact automorphism group and the unit ball in
$\Bbb C^n$ in 1977. Klembeck presented in 1978 a completely different proof
which is much more differential geometric. Then in 1979, Rosay strengthened the theorem to
a bounded domain such that there exists a sequence of automorphisms $\{f_j\}$ sending a
point $p$ to a smooth strictly pseudoconvex boundary point. Efimov  [4], using
Pichuk's rescaling method along with certain estimate on the Kobayashi metric, has extended
the above theorem of Rosay to unbounded domains in
$\Bbb C^n$. In an attempt to solve the rigidity for  Grauert tubes, the author and Ma
[10] also proved an extended version of the Wong-Rosay theorem to certain complete
hyperbolic  Stein manifolds. In this section
we will drop out the requirement on the complete hyperbolicity and  prove a version
of Wong-Rosay theorem for general domains inside complex manifolds.

\proclaim{Theorem 3.1}
Let $M$ be a domain in a  complex manifold $\hat M$. Suppose that there exist a point $p\in
M$ and a sequence of automorphisms $\{f_j\}\subset Aut(M)$ such that $f_j(p)\to q\in \pd M$, a
$C^2$-smooth strictly pseudoconvex point. Then $M$ is biholomorphic to the unit ball.
\endproclaim
\demo{Proof}

Let $U$
be a connected neighborhood of $q$ in 
$\hat M$ such that there exists a
biholomorphic map $\Phi$ from $U$ to a neighborhood $\Phi(U)$ of 0 in $\Bbb
C^n$ which maps  $q$ to the origin and such that the Kobayashi pseudodistance 
$d^{\sssize K}_{\sssize M}$ is a distance in $U$. Since $q$ is a
$C^2$-smooth strictly pseudoconvex point, we 
could choose suitable $\Phi$ and $U$ so that $\Phi(U)=E$ is a strictly convex domain in 
$\Bbb C^n$ and the domain
$$D=\Phi(U\cap M)=\{z\in E: \phi(z)<0)\}\subset \Bbb C^n$$ is strictly pseudoconvex with
the  defining function 
$$\phi(z)= \hbox{Re}\,z_n+|'z|^2+ \hbox{Re}\sum_{j=1}^nb_jz_j\bar z_n +o(|z|^2)$$
where $z=('z,z_n)$ and $b_j\in\Bbb C$.
Let 
$$ p_j=f_j(p),\ \ \ 2 r_j=d^{\sssize K}_{\sssize M}( p_j, M\backslash U).$$
From Lemma 2.2, we see that  $r_j\to \infty$ as $j\to \infty$. Let 
$$B_{\sssize K}(p,r)=\{w\in M : d^{\sssize K}_{\sssize M}(p,w)<r\}$$ denote the Kobayashi
ball of radius $r$ centered at $p$. Then 
$$f_j(B_{\sssize K}(p, r_j))=B_{\sssize K}(p_j, r_j)\subset U\cap M.$$
Hence,
$$K_j=\overline {B_{\sssize K}(p, r_j)}\subset f_j^{-1}(U\cap M).\tag3.1$$
For any $z\in M$, there exists a neighborhood $B_{\sssize K}(p, r_j)$ of $z$ and a $f_j\in Aut
(M)$ sending $B_{\sssize K}(p, r_j)$ to $U\cap M$. It follows from the biholomorphic invariance 
of
$F_{\sssize M}$ and (2.3) that 
$$F_{\sssize M}(z,\xi)= F_{\sssize M}(f_j(z),f_{j*}\xi)\ge
c\|\xi\|\;\;\text { for all }  z\in B_{\sssize K}(p, r_j),\;\xi\in
\Bbb C^n.$$
Hence, every point of $M$ is a hyperbolic point, i.e., $M$ is in fact a hyperbolic manifold.

Let $\zeta_j$ denote the unique point on $L=\Phi(\partial M\cap U)$ closest
to $a_j=\Phi(p_j)$ in the Euclidean distance.  Let $g_j$ denote the
composition of the 
translation that maps $\zeta_j$ to 0 and a unitary transformation,
taking the tangent plane $T_{\zeta_j}L$ to the plane
$\{\hbox{Re } z_n=0\}$. The mapping $g_j$ sends $D$ biholomorphically onto a bounded
domain $G_j$ with the defining function given by
$$\aligned 
\phi_j&=\phi\circ g_j^{-1} (z)\\
&= c_j\text{Re }z_n+d_j|'z|^2+  \text{Re }\sum_{k=1}^{n-1}b_{k}\mu^j_{k}\bar c_j
z_{k}\bar z_n+  \text{Re }b_n|c_j|^2|z_n|^2+0(|z|^2)\endaligned$$ where
$$c_j\to 1,\quad d_j\to 1,\quad \mu^j_{k}\to 1,\text { as } j\to\infty.$$
Since $a_j$ lies on the normal to $L$ at $\zeta_j$, the point 
$$g_j(a_j)=('0,-\delta_j).$$  By construction, $a_j\to 0$ and hence, $\zeta_j\to 0$ and
$\delta_j\to 0$ as
$j\to 0$. Let $F_j$ be the scaling function defined by 
$$F_j('z, z_n)=('z/\sqrt{\delta_j},z_n/\delta_j),$$ 
and let $h_j=F_j\circ g_j$. Then
$h_j(a_j)=('0,-1)$ and $h_j$ maps $D$ biholomorphically to a
domain $H_j\subset\Bbb C^n$ with the defining function
$$\aligned
\psi_j (w)&=(1/\delta_j)\phi\circ h_j^{-1} (w)\\
&=(1/\delta_j)\phi\circ g_j^{-1}\circ F_j^{-1} (w)\\
&=(1/\delta_j)\phi_{j}(\sqrt\delta_{j}\ 'w,\delta_j w_n)\\
&= c_j\text{Re }w_n+d_j|'w|^2+  \text{Re
}\delta_j^{\frac 12}\sum_{k=1}^{n-1}b_{k}\mu^j_{k}\bar c_j w_{k}\bar w_n\\
&\qquad+  \text{Re
}\delta_j \ b_n|c_j|^2|w_n|^2+\delta_j ^ {\varepsilon}0(|w|^2),\text{ for some
}\varepsilon>0.\endaligned$$
 The defining
functions $\psi_j$ converge, uniformly on compact subsets of $H_j$, to\
$\psi(w):= \hbox{Re 
}w_n+|'w|^2$ which is the defining function of the Siegel half-space $\Cal U$,
biholomorphically equivalent to the ball.   Let 
$$R_\nu= \{w\in\Bbb C^n: \hbox{Re } w_n+(1-\nu)|'w|^2<0\}.$$
These  domains satisfy  
$R_{\nu_1}\subset R_{\nu_2}$ for $\nu_1< \nu_2$.
Let
$$
\nu_j=\inf\{\nu>0: H_j\subset R_\nu \}.$$ 
Then $\nu_j\to 0\;\;\hbox{as}\;\;j\to\infty.$
Without loss of
generality we assume that $\nu_j<1/2$. 

Consider the maps 
$$\Psi_j=h_j\circ \Phi\circ f_j: \;\;
f_j^{-1}(U\cap M) 
\to H_j\subset R_{\nu_j}.$$  
By (3.1), each compact subset of $M$ is contained
in $f_j^{-1}(U\cap M) $ for sufficiently large $j$.
Since the domains $R_{\nu_j}$ are contained in $R_{1/2}$, which is
biholomorphic to the unit ball in $\Bbb C^n$. By Montel's theorem, there exists some
subsequence
$\{\Psi_{j_\nu}\}$ converges in the compact open topology to a holomorphic map $\Psi: M \to
R_{1/2}$. As $\nu_j\to 0$ we conclude that $\Psi(M)\subset R_\nu$ for
each $0<\nu<\frac 12$, hence  
 $\Psi:M\to \Cal U$  is  a holomorphic map. 

We would like to show that $\Psi$ is actually
a biholomorphism. Following the idea of Efimov [4], we first show that $\Psi$ is a locally
one-one map.

Since $\Psi_j (p)=('0,-1)$ for every $j$ and $H_j\to \Cal U$ as $j\to\infty$, we may assume
that there exists  a positive
$\beta<\frac 12$ such that the Euclidean ball $\hat B:=B(('0,-1),\beta)$ is contained in all
$\Psi_j(K_j)\cap\Cal U$ for $j$ sufficiently large. Since
$M$ is hyperbolic, $\overset {\sssize\circ }\to K_j $, the interior of
$K_j$ defined at (3.1), is hyperbolic as well. By the hyperbolicity of $M$,  the
holomorphic decreasing and the invariance property of a hyperbolic metric, we have 
$$\aligned
c_1\|\xi\|&\le F_{\sssize M} (p,\xi)\le F_{\overset {\sssize\circ }\to K_j} (p,\xi)\\
&= F_{\sssize  \Psi_j (\overset {\sssize\circ }\to K_j)} (('0,-1), d\Psi_j|_p(\xi))\\
&\le F_{\sssize \hat B } (('0,-1), d\Psi_j|_p(\xi))\\
&\le c_2 \|d\Psi_j|_p(\xi)\|;
\endaligned$$ the last inequality comes from (2.1). 
Hence there are constants $c_1$ and $c_2$ independent of $j$ such that 
$$\|d\Psi_j|_p(\xi)\|\ge
\frac {c_1}{c_2}\|\xi\|.\tag 3.2 $$ Since $\Psi_j$ are  biholomorphisms, $d\Psi_j$ are
nowhere zero for every $j$. The Hurwitz's theorem implies that the limit function $d\Psi$
is either nowhere zero or identically zero. The uniform estimate of (3.2) gives that 
$d\Psi|_p\ne 0$. Hence $d\Psi$ is nowhere zero and $\Psi$ is locally one-one. 

We claim
that it is in fact globally one-one. Suppose  there exist $x,y\in M$ such that $\Psi
(x)=\Psi (y)=s\in \Cal U$. Then there is  a Euclidean  ball $\Cal V=B(s,l)\subset \Cal U$ 
and two disjoint  neighborhoods $V_x$ and $V_y$ of $x$ and $y$ in
$M$ respectively  such that both $ \overline V_x$ and $\overline V_y$ are biholomorphic to
$ \overline{\Cal V}$ through $\Psi$. 
Choose $j$ sufficiently large so that $|\Psi_j-\Psi|<\frac {\ell}2$ in $V_x\cup V_y$.
  The boundary of the connected domain $\Psi_j (V_x)$ 
is less than $\frac {l}2$ away from the boundary of  
$\Cal V$ and $|\Psi_j (x) -s|< \frac {l}2.$
 Thus $s\in  \Psi_j (V_x).$
Similarly, $s\in  \Psi_j (V_y).$ That is $\Psi_j (V_x)\cap
\Psi_j (V_y)\ne\emptyset$ which contradicts to the injectivity of $\Psi_j$.

 Hence $\Psi$ is
a one-one holomorphic map, $M$ is biholomorphic to its image $\Psi
(M)\subset \Cal U$. As $\Cal U$ is   biholomorphic to the unit ball in $\Bbb C^n$, we view
 $M\approxeq \Psi(M)$ as a bounded domain in $\Bbb C^n$ where the Montel's theorem applies.

Consider the map
$$\Psi^{-1}_j: H_j\to f_j^{-1}(U\cap M)\subset M.$$ 
There exists a subsequence  $\Psi^{-1}_{\nu_j}$ converges to a holomorphic map 
$$\hat\Psi: \Cal U\to  M.$$ From the construction, it is clear that $\hat\Psi$ is the
inverse of $\Psi$. Hence $\Psi$ is a biholomorphism from $M$ to the ball $\Cal U.$\qed
\enddemo

\ 

\

\subheading { 4. General property of  $Aut(T^rX)$}

 Let $(X,g)$ be a real-analytic Riemannian manifold.  The {\it \gt} 
$$T^rX=\{(x,v):x\in X, v\in T_xX, |v|<r\}$$ is the collection of
tangent vectors of length less than $r$   equipped with the unique complex structure,
the {\it adapted complex structure}, which turns each leaves of the Riemannian foliation on
$T^rX\backslash X$ into holomorphic curves. That is, the map $f(\mu +i\tau)=\tau
\gamma'(\mu),\; \mu +i\tau\in \Bbb C$, is holomorphic, whenever it is defined, for any
geodesic $\gamma$ of $X$. There is a natural anti-holomorphic involution $\sigma$ fixing
every point of the center,
$$\sigma: T^rX\to  T^rX,\;\; (x,v)\to  (x,-v).$$
The length square function defined as $\rho (x,v)=|v|^2$ is \spsh\  and satisfies
the complex homogeneous Monge-Amp\`ere equation $(d d^c\rho)^n=0$ on $T^rX\backslash X$. The initial
condition for $\rho$ is that $\rho_{i\bar j}|_{\sssize X}=\frac 12 g_{i j}.$  We usually called the set
$(T^rX,X,g,\rho)$ a {\it Monge-Amp\`ere model }; $X$ the  center and $r$ the radius. 

For each real-analytic Riemannian manifold, there exists a  
$r_{max}(X)\ge 0$, maximal radius, such that the adapted complex structure is well-defined  
on
$T^{r_{max}}X$ whereas it blows up on   $T^sX$ for any $ s>r_{max}$. When the center
manifold $X$ is compact, it is clear that $r_{max}(X)>0$ since we could always paste the
locally defined  adapted complex structure together. For non-compact $X$, it is very likely
that the maximal possible radius degenerates to zero. However, when $(X,g)$ is  
homogeneous  the corresponding  $r_{max}(X)$ is always positive. The same holds for
co-compact $X$. The K\"ahler manifold $T^rX$ is a submanifold of $ T^{r_{max}}X$ for any
$r<r_{max}(X)$. 

It is clear from the construction of \gt s that the differential of an isometry of $(X,g)$
is actually an automorphism of $T^rX$, i.e., $Isom (X)\subset Aut
(T^rX)$. We ask whether  the converse holds or not. 

If $Isom (X)= Aut
(T^rX)$, i.e., every automorphism comes from the differential of an isometry of $X$, 
we say the \gt\ $T^rX$ has the {\it rigidity} property. In the compact case, this is
equivalent to the uniqueness of center: there is no other Riemannian manifold $(Y,h)$ such
that the complex manifold $T^rX$ is represented as $T^rY$, a \gt\ over $Y$.

When the center $X$ is compact, the rigidity of \gt s has been completely solved by
Burns-Hind in [2]. They prove that $Aut (T^rX)= Isom (X)$ for any compact 
real-analytic Riemannian manifold $X$ of any $r\le r_{max}(X)$. For the non-compact
centers, the only result obtained so far is contained in [8] [9] and [10]. The authors prove
the rigidity statement for
\gt s
$T^rX$ over any symmetric space $X$ of $r<r_{max}(X)$ provided that $T^rX$ is not
covered by the ball. 

In the compact case, one of the key observation toward the rigidity of \gt s is the
following: an automorphism
$f$ of the Grauert tube
$T^rX$ comes from the differential of an isometry   of the Riemannian manifold $(X,g)$ if
and only if
$f$ keeps the center $X$ invariant, i.e., $f(X)=X$. We will show in the following theorem
that Grauert tubes over non-compact centers share the same property. 
\proclaim{Theorem 4.1}
Let $f\in Aut (T^rX),\; r\le r_{max}$. Then $f=du$ for some $u\in Isom (X)$ if and only if
$f(X)=X$.
\endproclaim
\demo{Proof} One direction  is clear. Suppose $f(X)=X$.
We denote the Monge-Amp\`ere model over  $(X,g)$ as
$(\Omega,X,g,\rho)$ which means  $\Omega=T^rX$ is the Grauert tube over $X$ of radius $r$; 
$\rho$ is a non-negative
\spsh function satisfying the homogeneous complex Monge-Amp\`ere equation on $\Omega\backslash
X$; the zero set of $\rho$ is exactly $X$ and $(g_{i j})=2(\rho_{i\bar j}|_{\sssize X})$ is the
Riemannian metric on
$X$ induced from the K\"ahler form  ${i}\pd\bar\pd\rho.$  The automorphism
$f$ maps the Monge-Amp\`ere model
$(\Omega,X,g,\rho)$ onto the Monge-Amp\`ere model $(\Omega,X,k,\rho\cdot f^{-1})$ where $k$
is the Riemannian metric on $X$ induced from the K\"ahler form $
{i}\pd\bar\pd(\rho\cdot f^{-1})$. That is to say that $\Omega$ is a Grauert tube of radius
$r$ over the Riemannian manifold
$(X,k)$ as well. Let 
$$u=f|_{X}:(X,g)\to (X,k)$$
be the restriction of  $f$ to $X$. The map $u$ is an isometry from $(X,g)$ to $(X,k)$.
Actually, it is  an isometry of
$(X,g)$ since the Grauert tube $\Omega=T^r(X,g)=T^r(X,k)$ is the collection of tangent
vectors of
$X$ of length less than $r$ under both metrics $g$ and $k$. This forces the  metric $k$ to
be equal to the  metric
$g$.

Hence $u\in Isom (X,g)$ and $du\in Aut(\Omega).$ The automorphism $du\cdot f^{-1}$ of
$\Omega$  is the identity on a  maximal totally real submanifold $X$. Therefore
$$f=du$$
 on the whole $\Omega$.\qed
\enddemo
It is in general easier to deal with simply-connected manifolds. We will assume all of the
  manifolds $X$ are simply-connected for  proofs in Section 7. Before making
this assumption, we claim that it won't do any hurt for the general situations. Notice that
for any given $X$, we could always lift it to its universal covering $\tilde X$ which is
simply-connected. Let's  denote the \gt\ $T^rX$ as $\Omega$. Then the \gt\ $T^r\tilde
X=\tilde\Omega$ is the universal covering of $\Omega$. We have the following relation on
the rigidity of $\tilde\Omega$ and $\Omega$.
\proclaim{Lemma 4.2}
Let $X, \tilde X, \Omega$ and $\tilde\Omega$ be as above. Suppose there is a unique \gt\
representation for $\tilde\Omega$, then there is a unique \gt\
representation for $\Omega$
\endproclaim
\demo{Proof}If $\Omega=T^rX=T^rY$. Then $\tilde \Omega=T^r\tilde
X=T^r\tilde Y$ has two
\gt\ representations. Thus $\tilde X=\tilde Y$. 
On the other hand, $X=\tilde X/\Gamma$ for
some
 $\Gamma\subset Isom (\tilde X)$. Hence $\Omega=\tilde\Omega/\Gamma$. It follows that $Y=\tilde
Y/\Gamma=\tilde X/\Gamma=X$.\qed

\enddemo

Notice that the universal covering of a homogeneous Riemannian manifold is
homogeneous (c.f. Theorem 2.4.17 [17]). Hence the rigidity of \gt s over homogeneous
Riemannian  spaces  automatically holds once we prove it for simply-connected homogeneous
Riemannian  spaces.

\

\
\subheading {5. The complete hyperbolicity  of $T^rX$}

   It was proved by Sibony on Theorem 3 of [16] that a complex manifold admitting a bounded
\spsh function is hyperbolic. This implies that every \gt\ of finite radius is a hyperbolic
manifold. When the center  is compact and the radius $r<r_{max}$, the \gt\ $T^rX$ is a strictly
pseudoconvex domain sitting  inside the Stein manifold  $T^{r_{max}}X$ and hence  is
complete hyperbolic. There is no such generality for   non-compact cases. So far, the only
known characterization is over  co-compact manifolds.

When $X$ is co-compact, i.e.,
$ X/\Gamma$ is compact for some
$\Gamma\subset Isom(X)$, the \gt\ $T^rX$ is   simply the universal covering of $T^r(X/\Gamma)$
and hence is completely hyperbolic in case $r<r_{max}$. We will show in this section that  
a \gt\  constructed over a homogeneous Riemannian  manifold of radius less than the maximal
 is complete hyperbolic. 

 A sequence $f_{j}\in Hol(N,W)$ is said to be {\it compactly divergent } if for
each compact subset
$K\subset N $ and each compact $K'\subset W$ there is a $j_0=j_0(K,
K')$ such that the set
$f_{j} (K)\cap K'$ is empty for each $j\ge j_0$. A complex manifold $W$
is said to be {\it taut}  if each sequence in $Hol(N,W)$ contains
a subsequence that either converges to an element of $Hol(N,W)$ in
the compact open topology, or diverges compactly.   It is known that a complete
hyperbolic manifold is taut and a taut manifold is hyperbolic.

A \Rm\ is said to be {\it homogeneous} if its isometry group acts transitively on it. It is a
classical result that a \Rm\ $(X,g)$ is complete if and only if every bounded subset
is relatively compact in $X.$ Let $d$ denote the distance for the metric $g$ in $X$ and
$B(x,r)=\{y\in X: d(x,y)<r\}$ denote the  ball of radius $r$
centered at $x\in X$.
\proclaim{Lemma 5.1}
Let $X$ be a non-compact homogeneous \Rm, $x\in X$. Then for any $R>0$ there exist $y\in
X$ such that $d(x,y)>R.$
\endproclaim
\demo{Proof}
Suppose there exist a $x\in X$ and a  $R>0$ such that $d(x,y)<R$ for all $y\in X$. Then
$X\subset B(x,R)$. Since $B(x,R)$ is  a bounded subset of $X$; it is relatively
compact in
$X$. Thus, $\overline {B(x,R)}$ is a compact subset of $X$ and 
$$\overline {B(x,R)}\subset X\subset B(x,R)$$
which implies that $X=\overline {B(x,R)}=B(x,R)$ is compact, a contradiction.\qed
\enddemo
Let $\Omega=T^rX$ be a \gt\ of $r<r_{max}$ and let
 $d_{\sssize K}$ denote the Kobayashi distance on $\Omega$. The boundary of $\Omega$
contains two parts. The first part is consisted of smooth strictly pseudoconvex points
$\{(x,v):\; v\in T_xX,\; |v|=r\}$. The second part comes from the boundary of the
Riemannian manifold $X$; it is the set $\{(x,v):\; v\in T_xX,\; |v|<r, x\in B^{\sssize C}(p,R),
 p\in X, \forall R>0\}$.

From Lemma 2.2 we
know that for any $p\in\Omega$, $d_{\sssize K}(p,q)\to\infty$ as $q$ approaches  a smooth
strictly pseudoconvex boundary point. A metric is complete means that every bounded subset is relatively
compact. To show the complete hyperbolicity of
$\Omega$, it is therefore sufficient to show that $d(p,y)\to\infty$ as $y$ goes to the boundary of $X$.
\proclaim{Theorem 5.2}
Let $X$ be a  homogeneous \Rm, $\Omega=T^rX$ be a \gt\ of $r<r_{max}$. Then
$\Omega$ is complete hyperbolic.
\endproclaim
\demo{Proof}
The statement  is clear for compact $X$ by the existence of a \spsh exhaustion
function.  We shall prove the non-compact case. 

For $p\in X$, let $S(p,n)=\{x\in X: d(p,x)=n\}$ be the 
$n$-sphere in $X$ around $p$. By Lemma 5.1, $S(p,n)\ne\emptyset$ for all $n\in\Bbb N$. Picking a
point $p_2\in S(p,2)$, there exists a unique minimal geodesic $\gamma_2(t)$, parametrized by the
arclength, joining $p$ to $p_2$. The geodesic $\gamma_2(t)$ also serves as the minimal geodesic
for any two points $\gamma_2(a)$ and $\gamma_2(b)$, $0\le a\le b\le 2$. Then 
$p_1:\; =\gamma_2(1)\in S(p,1)$ and $d(p_1,p_2)=1$.  Let 
$$\mu_1=\min\{d_{\sssize K} (p,x):x \in S(p,1)\}>0.$$  Since $p_1\in S(p,1)$, 
 it follows that
$$d_ {\sssize K} (p,p_1)\ge\mu_1.\tag 5.1$$
 Since $X$ is homogeneous, there exist $g_j\in Isom (X)$ such that 
$$ g_j(p_j)=p,\quad j=1,2.$$
Hence
$$d(p,g_1(p_{2}))=d(g_1(p_1),g_1(p_{2}))=d(p_1,p_{2})=1.$$
This implies that $g_1(p_{2})\in S(p,1)$ and therefore
$$d_ {\sssize K} (p,g_1(p_{2}))\ge\mu_1.\tag 5.2$$
As the Kobayashi distance is biholomorphically
invariant and $g_j\in Isom (X)\subset Aut(\Omega)$, we have
 $$ d_ {\sssize K} (p_1,p_{2})=d_{\sssize K}
(g_1(p_1),g_1(p_{2}))
=d_{\sssize K} (p,g_1(p_{2}))\ge\mu_1.\tag 5.3$$
 Let $B_{\sssize K}(x,\mu_1)$ denote the restriction of the Kobayashi ball of radius
$\mu_1$ to $X$, i.e.,
$$B_{\sssize K}(x,\mu_1)=\{ y\in X: d_{\sssize K} (x,y)<\mu_1\}.$$
Then $B_{\sssize K}(p,\mu_1)\subset B(p,1)$. Therefore  we have, for $j=1,2,$
$$\aligned
B_{\sssize K}(p_j,\mu_1)&=B_{\sssize K}(g_j^{-1}(p),\mu_1) 
 =g_j^{-1}B_{\sssize K}(p,\mu_1)\\&\subset g_j^{-1}B(p,1)
=B(g_j^{-1}(p),1)=B(p_j,1).\endaligned\tag 5.4$$
For $x\in B(p,1)$, $d(p_{2},x)\ge d(p_{2},p)-d(p,x)\ge 1.$
Hence  $$B(p ,1)\cap
B(p_{2 },1)=\emptyset.\tag 5.5$$
(5.4) and (5.5) immediately imply that 
 $$ B_{\sssize K}(p,\mu_1)\cap B_{\sssize K}(p_{2},\mu_1)=\emptyset.\tag 5.6$$
Thus, $d_{\sssize K}(p,p_{2})\ge 2\mu_1.$ Since $p_2$ could be taken from   arbitrary points in
$S(p,2)$, 
 we have 
$$d_{\sssize K}(p,q)\ge 2\mu_1\;\;\text{ for all } q\in S(p,2).\tag 5.7$$
Repeating the same precess by taking   $q_4\in S(p,4),$ the corresponding minimal geodesic
 joining $p$ and $q_4$ be $\gamma_4$, $q_2=\gamma_4(2)\in S(p,2)$ and
$$\mu_2=\min\{d_{\sssize K} (p,x):x \in S(p,2)\}.$$ By (5.7), $\mu_2\ge 2\mu_1$.
We get $$d_{\sssize
K}(p,q)\ge  2 \mu_2\ge  2^2\mu_1,\;\;\text{ for all } q\in S(p,2^2).\tag 5.8$$
By inductive argument, we conclude that 
$$d_{\sssize
K}(p,q)\ge  2^n\mu_1,\;\;\text{ for all } q\in S(p,2^n).\tag 5.9$$
 Finally, we claim that $\Omega$ is complete hyperbolic, i.e., every bounded subset is relatively
compact in
$\Omega$ or equivalently,
$$\lim_{q\to \pd\Omega}d_{\sssize K}(z,q)=\infty.\tag 5.10$$
It was already shown on Lemma 2.2 that if $q$ goes to a smooth strictly pseudoconvex boundary
point then (5.10) holds, so it is sufficient to check those $q\in X$ staying away from any bounded set
in $X$. If  $q$ stays outside of $B(p,2^n)$  we  conclude that $d_{\sssize
K}(p,q)\to\infty$ as  $n\to\infty$. Moreover, for any
$z\in
\Omega$,
$$d_{\sssize K}(q,z)\ge d_{\sssize K}(q,p)-d_{\sssize K}(p,z)\to\infty$$
as well.  Hence the Kobayashi metric $d_{\sssize K}$ is complete.\qed

\enddemo

\

\subheading { 6.  The equality of $Aut_0(T^rX)=Isom_0 (X)$} 

Since a  \gt\
of finite radius  is hyperbolic, its
automorphism group is a Lie group. 
It is clear
from the construction of \gt s that the differential $du$ is an
automorphism of $T^rX$ for any $u\in Isom (X)$. Hence
the isometry group $Isom(X)$ is a Lie subgroup of  $ Aut
(T^rX)$.  Without any ambiguity, we also use the same
symbol $u$ to represent its differential $du$ in $Aut (T^rX)$.
 
From now on we shall assume that $X$ is a real-analytic Riemannian manifold such that the
adapted complex structure exists on 
$T^rX$. 
    Let
$$T^r_pX=\{(p,v): v\in T_pX, \;|v|<r\}$$ be the fiber through the point $p\in X$ and  $\sigma$
be the  anti-holomorphic involution 
$$\sigma: T^rX\to  T^rX,\;\; (p,v)\to  (p,-v).$$ 
When $r<r_{max}$, the boundary $\partial T^r_pX$ is consisted of smooth strictly pseudoconvex
points since it is  locally defined by
$\{\rho^{-1}(r^2)\}$.

 Let $I$ denote the
isometry group of $X$ and 
$G$ denote the automorphism group of $T^rX$. Then the group 
$$\hat G:=  G\cup \sigma \cdot  G$$ is again a Lie group acting on $T^rX$. Let $\hat\Cal G$ be
the Lie algebra of $\hat G$. Each  $\xi\in\hat\Cal G$ could be viewed as a vector
field in $T^rX$; $\xi (z)=\frac {d}{dt}|_{t=0}(exp\ t\xi)(z)$.

\proclaim {Lemma 6.1}
For each $\xi\in \hat \Cal G$, there is an $\eta\in \hat \Cal G$ corresponding to
$\xi$ such that for any $t\in\Bbb R$, $\exp t\eta: T^r_pX\to T^r_pX,\ \forall p\in X$.
\endproclaim
\demo{Proof}
Define $\eta\in \hat\Cal G$  as 
$$\eta=\xi-\frac {d}{dt}|_{t=0} (\sigma\cdot (exp\ t\xi)).$$
In local coordinates,  $\sigma (z)=\bar z$. Then
$$\eta (iy)= 2i\ \text{Im } \xi (iy).$$
It follows that $\eta$ is tangent to the fiber $T^r_pX$ and the result follows.\qed
\enddemo

Fix a point $p\in X$ and take $U$ to be a small neighborhood of $p$ in $X$. Take
 $$D=\{(x,v):x\in U, v\in T_xU,|v|<r\}$$ be a domain in $T^rX$ with the induced complex
structure from $T^rX$.
 The domain $D$  could be equipped  with the 
metric $d$  induced from the Kobayashi metric $d_{\sssize K}$ of $T^rX$,
i.e., define the metric $d(z,w):=d_{\sssize K}(z,w),\forall z,w\in D.$

 For a fixed $\eta\in\hat\Cal G$ as  in the statement of the above lemma,
we consider the restriction maps $exp\ (t\eta)|_{\sssize D}$
and set
$$ 
\Cal F=
 \{exp\ (t\eta)|_D; t\in\Bbb R\}\subset C(D,T^rX), $$
where $ C(D,T^rX)$ denotes the set of all  continuous maps from $D$
to $T^rX$.
\proclaim {Lemma 6.2} Suppose $T^rX$ is not biholomorphic to the ball and $r<r_{max}$, then
$\Cal F$ is a compact subset of $C(D,T^rX).$
\endproclaim
\demo{Proof}
It is clear that $\Cal F$ is closed in $C(D,T^rX)$.
By Lemma 6.1, $f(D)\subset D$ for all $f\in \Cal F.$
As $d_{\sssize K}$ is an
invariant metric and $exp\ (t\eta)\in Aut (T^rX)$, we have 
$$d(exp\ (t\eta)(z),exp\ (t\eta)(w))=d_{\sssize
K}(exp\ (t\eta)(z),exp\ (t\eta)(w))= d_{\sssize K}(z,w)= d(z,w)$$ for all $z,w\in D$ and $t\in
\Bbb R$. This shows that $\Cal F$ is equicontinuous. We then claim that for every
$z\in D$, the set  $\Cal F(z):=\{exp\ (t\eta)(z):t\in \Bbb R\}$ is
relatively compact in
$T^rX$. Suppose not, for some $z=(p,v)\in D$ the set of points $\{exp\ (t\eta) (z)\}$ approach to
the boundary of
$T^r_pX$ which is consisted of  smooth strictly pseudoconvex points. By Theorem 3.1, this forces
$T^rX$ to be the ball which is a contradiction. Therefore $\Cal F(z)$ is a relative
compact subset of $T^rX$. By the Ascoli theorem (c. f. [19]), $\Cal F$ is  compact
in
$C(D,T^rX)$.\qed
\enddemo

\proclaim {Lemma 6.3}
For each $\xi\in \hat \Cal G$, the vector field $\xi$ is tangent to $ X$.
\endproclaim
\demo{Proof}By lemmas 6.1 and 6.2, $\Cal F$ is a connected compact Lie group acting
symmetrically on each fibre $T_p^rX$. The action is symmetric as defined on $\S 4$ of [10]
since if $g\in \Cal F$ then $\sigma\cdot g\cdot \sigma (y)=-g(-y)$. Applying Theorem 4.1 of
[10], every   $p\in X$ is a fixed point of any $g\in \Cal F$. From the construction of
$\eta$ in Lemma 6.1, this implies that $$ \text{Im } \xi (p)=-\frac {i}2\eta
(p)=0,\;\;\forall p\in X.$$
Therefore $\xi$ is tangent to $X$.
\qed
\enddemo 
Denote  the identity component of $Isom (X)$ by $I_0$ and 
the identity component of $Aut(T^rX)$ by $G_0$. We prove the following for Grauert tubes over any
real-analytic Riemannian manifold.
\proclaim{Theorem 6.4}
Let $X$ be a real-analytic Riemannian manifold such that $r_{max}(X)>0$.  Then
$Aut_0 (T^rX)=Isom_0 (X)$ for any $r<r_{max}$ provided that $T^rX$ is not
covered by the ball.
\endproclaim 

\demo{Proof}
There exist a neighborhood $V$ of 0 in $\hat \Cal G$ and a neighborhood $U$ of $id$ in $\hat
G_0$ such that the exponential map is a diffeomorphism from $V$ to $U$. That is, for every
$f\in U$, there exists an $\xi\in V$ such that $f=exp\ \xi$. By Lemma 6.3, $f(X)\subset X.$
As $X$ is a closed submanifold and $f$ is a diffeomorphism of $T^rX$, $f(X)$ is then a
closed submanifold of the connected manifold $X$ of the same dimension, which implies that
$f(X)=X.$ 

By Theorem 4.1,
$f\in I_0.$ Hence $U\subset I_0\subset \hat G_0$ and we conclude that the manifolds $I_0$ and
$\hat G_0$ have the same dimension. Since $I_0$ is a Lie subgroup of $\hat G_0$, $I_0$ is a
closed submanifold of $\hat G_0$, therefore
$I_0=\hat G_0$. On the other hand, it is clear from the definition of $\hat G$ that $\hat G_0=
G_0$. We conclude that
$I_0=G_0$. 
\qed
\enddemo

\

\subheading { 7.  The rigidity of \gt s over homogeneous centers} 

For compact $X$, Burns-Hind [2] have proved that the isometry group   $Isom (X)$ of $X$   is
isomorphic to the automorphism group $ Aut (T^rX)$ of the \gt\   for any radius $r\le r_{max}$.
Burns [1] also proved that this isomorphism is equivalent to the uniqueness of the center.

For the non-compact cases, the only known rigid \gt s are \gt s 
constructed over locally symmetric spaces of $r< r_{max}$ in [10] provided that  the \gt s are not
covered by the ball. In [9], the author further proved  the uniqueness of the center also holds for
the above non-compact cases.

 We will prove in this section that not only the rigidity but also the uniqueness holds for 
\gt s over non-compact homogeneous \Rm s  of $r< r_{max}$; the only exception occurs when the
\gt\
 is biholomorphic to the ball. It was proved by the author in [8] that  $T^rX$ is biholomorphic
to
$B^n\subset\Bbb C^n$ if and only if  $X$ is the real hyperbolic space $\Bbb H^n$ of curvature $-1$
and the radius
$r=\frac {\pi}{4}$. (The Monge-Amp\`ere solution $\rho$ in this article is half of the $\rho$ in
[8], hence the radius changes.)
Apparently, the automorphism group of $B^n$ is much larger
than the isometry group of $\Bbb H^n$.

We assume in this section
that $X$ is a simply-connected Riemannian homogeneous space and $\Omega=T^rX$ is the Grauert tube
over
$X$ of radius $r$. Recall that a Riemannian manifold  $X$ is    homogeneous  if the isometry
group acts transitively on $X$. Denote $I=Isom (X)$ and $G=Aut(T^rX)$; $I_0$ and $G_0$ the
corresponding identity components. Notice that the homogeneity immediately implies that the
$X=I_0(p)$ for any
$p\in X$. 

For any given $f\in G$, the set $Y=f(X)$ equipped with the push-forward metric coming from $X$
is a homogeneous center of the Grauert tube $\Omega$. We claim that  $Y$ crosses through every
fiber $T_p^rX$.  
\proclaim{Lemma 7.1}
$Y\cap T_p^rX\ne\emptyset$ for any $p\in X$.
\endproclaim
\demo{Proof}
Since $X$ is a connected homogeneous space, the fact that
$I_0=G_0$ is a normal subgroup of $G$ would imply that for any $p\in X$, we have
$$Y= f(X)=f(I_0\cdot p)=I_0\cdot f(p)\tag 7.1 $$
is an $n$-dimensional $I_0$-orbit as well. Let 
$$ \pi: \Omega\to  X, \;\; \pi (p,v)=p$$ be the projection map. It is clear that $\pi$ is
$I_0$-equivariant since for any $(x,v)\in \Omega$ and $g\in I_0$, we have
$$\pi (g\cdot (x,v))=\pi (g\cdot x, g_* v)=g\cdot x=g\cdot \pi (x,v).$$
By (7.1),
$$\pi (Y)=\pi (I_0\cdot f(p))=I_0\cdot \pi(f(p))=X. $$ Therefore $Y$ intersects every fiber
$T_p^rX$. The lemma is proved.\qed

\enddemo
The Lie group $I $ is a  subgroup of $G $. We consider the quotient space
$G /I $ and examine the index of this coset space. 

\proclaim{Proposition 7.2}
The index of $I $ in $G $ is  finite.
\endproclaim
\demo{Proof}
Let $\{g_j I \}, g_j\in G $, be a sequence in the coset space $G /I $. By Lemma 7.1,
$g_j(X)$  has non-empty intersection with the fiber $T^r_pX$. Take a
point $$q_j\in g_j(X)\cap  T^r_pX.$$ Since $g_j^{-1}(q_j)\in X$ and $I $ acts
transitively on $X$, there exists a
$h_j\in I $ sending the point $p\in X$ to  $g_j^{-1}(q_j)\in X$. Let 
$$f_j=g_j\cdot h_j\in G $$ then $f_j(p)=q_j\in T^r_pX$ for all $j$. Since the boundary of
$T^r_pX$ is consisted of smooth strictly pseudoconvex points  and the \gt\ $\Omega$ is not
the ball, the generalized Wong-Rosay
theorem in $\S 3$ implies that  no subsequence of $\{f_j\}$ could diverge compactly. 
By Theorem 5.2, $\Omega$ is a taut manifold which asserts that there exists a
subsequence of
$\{f_j\}$ converges to some
$f \in G $ in the topology of $G $. Hence, some subsequence of $\{g_j I \}$ converges to
$g I $ in the topology of $G /I $. Thus, $G /I $ is compact. Since $I$ and $G$ have the same
identity components, 
$I$ is open in $G$. The compactness  implies that  the index of $G /I $ is finite.\qed
\enddemo
Suppose the Grauert tube $\Omega=T^rX$ could be represented as a Grauert tube of the same
radius of another center $Y$, i.e.,
$\Omega=T^rY$. We  use the
notations $(X,\sigma)$ and $(Y,\tau)$  to 
 indicate that the two anti-holomorphic involutions of the Grauert tube $\Omega$ with respect to
the centers
$X$ and $Y$ are
$\sigma$ and $\tau$ respectively. We call $(Y,\tau)$ a {\it homogeneous center} if $Y$ is a
homogeneous space as well. Regrouping the terms, we have the following relation:
$$\aligned 
&\tau\cdot (\sigma\cdot\tau)^{2l}=(\sigma\cdot\tau)^{-l}\cdot \tau\cdot
(\sigma\cdot\tau)^{l}\\
&\tau\cdot (\sigma\cdot\tau)^{2l-1}=(\tau\cdot(\sigma\cdot\tau)^{l-1})\cdot \sigma\cdot
(\tau\cdot(\sigma\cdot\tau)^{l-1})^{-1}
\endaligned\tag 7.2$$

\proclaim{Proposition 7.3}
Let  $(Y,\tau)$   be a homogeneous center of $\Omega$  and let $k$ be the least
positive integer  such that  $(\sigma\cdot\tau)^k\in I $. Then
$k$ is odd and $(\sigma\cdot\tau)^k=id$. 
\endproclaim
\demo{Proof}
Since the index of $I $ in $G $ is finite there exists the least integer $k$ with
$(\sigma\cdot\tau)^k\in I .$ Write $(\sigma\cdot\tau)^k=du$ for some $u\in I $. By (7.2),
$$\sigma\cdot du=\tau\cdot (\sigma\cdot\tau)^{k-1}$$ is an  anti-holomorphic involution in
$\Omega$ with fixed point set  $Z$ where 
$Z=\tau\cdot (\sigma\cdot\tau)^{\frac {k-2}2}(X)$ when $k$ is even; $Z=(\sigma\cdot\tau)^{-\frac
{k-1}2}(Y)$ when  $k$ is odd. 

$(Z,\; \sigma\cdot du)$ is a homogeneous center, with the pushed forward metric coming from
$X$ or $Y$, of the
\gt\
$\Omega$. Suppose there is
$z\in Z-X, \;z=(x,v)\in T_xX,v\ne 0.$ Since $z$ is a fixed point of $\sigma\cdot du$,
$$(x,v)=\sigma\cdot du (x,v)=\sigma\cdot (u(x),u_*v)=(u(x),-u_*v).$$
We have $u(x)=x,\; v=-u_*v$. Hence the whole interval $$L=\{(x,\alpha v)\in
T^r_{x}X:|\alpha v|<r\}$$ is fixed by
$\sigma\cdot du$ which implies that  $Z\supset L$. Since $Z$ is a homogeneous space, there
exist $\{f_j\}\in Isom (Z)\subset G$ such that $$f_j(x)=(x,\frac {jr}{(j+1)|v|}v),\;j\in
\Bbb N.$$ Then $f_j(x)\to (x,\frac {r}{|v|}v)$, a smooth strictly pseudoconvex boundary
point of
$\Omega$, as $j\to\infty$. By Theorem 3.1, $\Omega$ is then biholomorphic to the ball, a
contradiction. Therefore $Z=X,\;
\sigma\cdot du=\sigma$ and hence $du= (\sigma\cdot\tau)^k=id$.

Suppose $k=2l$ is even. By (7.2), $\sigma=(\sigma\cdot\tau)^{l}
\cdot\sigma\cdot (\sigma\cdot\tau)^{-l}.$
Thus, $$\sigma\cdot (\sigma\cdot\tau)^l=(\sigma\cdot\tau)^l\cdot\sigma.$$ The 
$n$--dimensional closed submanifold
$(\sigma\cdot\tau)^l(X)$ is the fixed point set of   the anti-holomorphic involution
$\sigma$. By the uniqueness of the fixed point set, 
$(\sigma\cdot\tau)^l(X)=X$. By Theorem 4.1, 
$(\sigma\cdot\tau)^l\in  I $ which is  a contradiction since $l<k$.\qed

\enddemo

We also prove a proposition similar to the Theorem 2bis of [1].
\proclaim{Proposition 7.4}
 For any homogeneous center
$(Y,\tau)$ of
$\Omega$, there exists an $f\in G $ such that $f(X)=Y$. Furthermore, there are at most a
finite number of homogeneous centers  in $\Omega$.
\endproclaim
\demo{Proof}
Let $k=2l+1$ be the odd integer as stated in Proposition 7.3 and let
$f:=(\sigma\cdot\tau)^l\in G
$. The condition
$(\sigma\cdot\tau)^k=id$ holds if and only if
$$\tau=(\sigma\cdot\tau)^{2l+1}\cdot \tau=(\sigma\cdot\tau)^l
\cdot\sigma\cdot (\sigma\cdot\tau)^{-l}=f\cdot\sigma\cdot f^{-1}.$$ Thus, $f(X)$ is the
fixed point set of $\tau$, i.e., 
$f(X)=Y$.

 Since $ G/I$ has finite index $h$, there exist $\{g_j\in G\}_{j=1}^h $ such that 
$G$ is the disjoint union of 
$g_j I, j=1,\cdots,h.$ Then for any $\eta\in G, \eta(X)\in \{g_j(X)\}_{j=1}^h$ where
$(g_j(X),\ g_j\cdot\sigma\cdot
g_j^{-1}) $ are homogeneous centers of $\Omega$. Suppose there exists a homogeneous center
$(W,\hat\tau)$
 other than
$\{(g_j(X),(g_j\cdot\sigma\cdot g_j^{-1})\}_{j=1}^h,$ then there is an $\hat f\in G$ such that
$\hat f (X)=W,$ a contradiction.
\qed 
\enddemo

\proclaim{Lemma 7.5}
 $Isom (X)  = Aut (\Omega) $ if and only if there is a unique homogeneous center $(X,\sigma)$ for
$\Omega$. 
\endproclaim
\demo{Proof}
Suppose there exists another homogeneous center $(Y,\tau)$, then we could find 
  $f\in G $ such that $Y=f(X)$  by Proposition 7.4. The condition $I = G $
implies that $Y=X$.

Suppose there exists $\zeta\in G \backslash  I $ sending $(X,\sigma)$ to another homogeneous
center
$(\zeta(X),\;\zeta\cdot\sigma\cdot \zeta^{-1})$.  
 By the assumption,
$\zeta(X)=X$. Thus, by Theorem 4.1,
$\zeta\in I $, a contradiction.\qed  
\enddemo
Recall that a \gt\ $\Omega=T^rX$ is said to have a unique center $X$ if this is the unique
\gt\ representation for $\Omega$, i.e., there is no other $Y$ such that $\Omega=T^rY$.
Finally we prove the main theorem of this paper, the rigidity and the uniqueness of the center.

 \proclaim{Theorem 7.6}
Let $X$ be (or be a quotient manifold of) a real-analytic homogeneous Riemannian manifold. 
 Let
$\Omega=T^rX$ be a
\gt\ of radius
$r<r_{max}$ such that
$\Omega$ is not covered by the ball.  Then
$Isom(X)= Aut(\Omega)$ and $\Omega$ has a unique center.
\endproclaim
\demo{Proof} Since the proof of Proposition 7.2 works through when replacing $I$ by $G_0$, it
is easy to see that the index of $G/G_0$ is finite.

Following the  same arguments in the proof of
Lemma 7.4 in [10]. We are able to construct an $\hat G$-invariant strictly plurisubharmonic
non-negative function $$\psi(z)=\sum_{j=1}^k \rho(g_j(z))$$ in $\Omega$ where $\{g_1,\dots, g_k\}\in G$
so that $G/G_0=\{g_jG_0: j=1,\dots, k\}$. As $G_0$ acts transitively on $X$, the tangent
space $T_z(T^rX)$ could be decomposed as, for any $z\in T^rX$,  $$T_z(T^rX)=T_z(G_0 \cdot
z)+T_z(T_{\pi(z)}^r X).$$ Since $\psi$ is constant in $G_0 \cdot z$, every critical point of
the function $f_z:= \psi|_{T_{\pi(z)}^r X}$ is a critical point of $\psi$ and every critical
point of $\psi$ occurs at the critical points of the functions $f_z$.

 As $\psi$ is strictly plurisubharmonic, the above
decomposition implies that the real Hessian of $f_z$ is positive definite on the tangent
space 
$T_z(T_{\pi(z)}^r X).$ Since $f_z$ is proper on the  fiber, it follows that there is
exactly one critical point of $f_z$ which turns out to the the minimal point. Since $\psi\cdot
\sigma=\psi$, the minimum of $f_z$ occurs at $\pi (z).$ That is to say that the set of critical
points of $\psi$ is exactly  $X$.

Suppose there is another homogeneous center $Y.$ By proposition 7.4, $Y=F(X)$ for some $F\in G$.
Since
$\psi$ is
$G$-invariant, every point of 
$Y$ is a  critical point  of $\psi$. Therefore, $Y\subset X$. We conclude that  $X=Y$ and
$Isom(X)= Aut(\Omega)$ is proved following from Lemma 7.5.

We then prove that there is exactly one center in $\Omega$. Suppose there is a  non-homogeneous
center $W$ of the Grauert tube $\Omega$. This $W$ has to be connected and $dim\ W=dim\ X$.

 Let $z=(x,v)\in W$. Since $Isom_0(W)=Aut_0 (\Omega)=Isom_0(X)$,
for any $g\in Isom_0(W)$ there exists $u\in Isom_0(X)$ such that 
$$g\cdot(x,v)=dg\cdot(x,v)= du\cdot(x,v)=(u(x), u_*v).$$
Let $g$ run through the Lie group $Isom_0(W)$, the corresponding $u$ will then go through $Isom_0(X),$ 
we obtain
$$W\supset Isom_0(W)\cdot (x,v)=Isom_0(X)\cdot (x,v).$$ 
By the homogeneity of $X$, the left hand side is a closed submanifold of dimension $\ge n$. Therefore,
$Isom_0(W)\cdot (x,v)$ is a closed submanifold of dimension $n$ in $W$. As $W$ is connected, we conclude
that $W= Isom_0(W)\cdot (x,v)$ and hence $W$ is homogeneous, a contradiction.

\qed
\enddemo

\remark {Remark}\roster
\item Since every symmetric space is homogeneous, 
results in this article cover those in [9] and [10].  
\item  It was shown in [3] that \gt s over non-compact symmetric spaces of rank one  
of  the maximal radius
are not rigid.
\endroster
\endremark

\enddocument